\documentclass[journal]{journal}
\usepackage{graphicx}%
\usepackage{subcaption}
\usepackage{multirow}%
\usepackage{amsmath,amssymb,amsfonts}%
\usepackage{amsthm}%
\usepackage{mathrsfs}%
\usepackage{xcolor}%
\usepackage{textcomp}%
\usepackage{manyfoot}%
\usepackage{booktabs}%
\usepackage{algorithm}%
\usepackage{algorithmicx}%
\usepackage{algpseudocode}%
\usepackage{listings}%
\usepackage[absolute,overlay]{textpos}

\theoremstyle{thmstyleone}%
\theoremstyle{thmstyletwo}%

\theoremstyle{thmstylethree}%
\ifCLASSINFOpdf
\else
\fi
\begin{document}
%
\title{From Initial Data to Boundary Layers: Neural Networks for Nonlinear Hyperbolic Conservation Laws}
%
%
%

\author{Igor CIRIL$^1$ \and Khalil HADDAOUI$^2$ \and Yohann TENDERO$^3$}

%
%

\markboth{Journal of \LaTeX\ Class Files,~Vol.~6, No.~1, January~2007}%
{Shell \MakeLowercase{\textit{et al.}}: Bare Demo of IEEEtran.cls for Journals}
%



\maketitle
\begin{textblock*}{15cm}(1.5cm,27cm)
\footnotesize
\begin{tabular}{@{}ll}
$^1$ & DRII, Institut Polytechnique des Sciences Avancées, Ivry-sur-Seine, France \\
$^2$ & CerebraQuant Solutions, Meudon, France \\
$^3$ & École Centrale d'Électronique, Lyon, France \\
\end{tabular}
\end{textblock*}

\thispagestyle{empty}

\begin{abstract}
We address the approximation of entropy solutions to initial-boundary value problems for nonlinear strictly hyperbolic conservation laws using neural networks. A general and systematic framework is introduced for the design of efficient and reliable learning algorithms, combining fast convergence during training with accurate predictions. The methodology that relies on solving a certain relaxed related problem is assessed through a series of one-dimensional scalar test cases. These numerical experiments demonstrate the potential of the methodology developed in this paper and its applicability to more complex industrial scenarios.
\end{abstract}

\begin{IEEEkeywords}
Physics-Informed Neural Networks, Conservation Laws, Boundary Value Problems
\end{IEEEkeywords}

%
\IEEEpeerreviewmaketitle

\section{Introduction}\label{intro}

Nonlinear hyperbolic conservation laws play a central role in the mathematical modeling of physical systems where transport and wave propagation phenomena dominate. These equations arise in a wide variety of applications, including compressible fluid dynamics, multiphase flows, shallow water equations, traffic and biological transport~\cite{leveque2002finite, smoller1994shock}. A defining feature of such equations is the spontaneous formation of discontinuities, or shock waves, in finite time—even when starting from smooth initial data. In this regime, classical solutions break down, and one must work in the framework of weak solutions. Among these, only entropy-admissible solutions are physically meaningful, satisfying additional inequalities that enforce consistency with the second law of thermodynamics~\cite{dafermos2016hyperbolic}.

Over the past decades, the numerical approximation of entropy solutions has been a major topic in computational mathematics. Classical methods, such as finite volume schemes~\cite{godlewski2021numerical, eymard2000finite, bouchut2004nonlinear}, discontinuous Galerkin methods~\cite{cockburn2001runge}, or front-tracking techniques, are specifically designed to capture shock dynamics while respecting entropy conditions and handling complex geometries. These schemes are mathematically grounded and efficient, but they typically rely on structured grids, high-resolution Riemann solvers, and ad hoc strategies to enforce physical boundary conditions and couple subsystems.

In parallel, recent years have witnessed a surge of interest in machine learning-based methods for solving partial differential equations (PDEs). In particular, Physics-Informed Neural Networks (PINNs)~\cite{raissi2019physics} offer a promising, mesh-free approach for approximating PDE solutions by minimizing a composite loss function that encodes the governing equations, initial data, and boundary constraints. While PINNs have shown strong performance for smooth solutions of elliptic and parabolic PDEs, their application to nonlinear hyperbolic problems remains limited. The main obstacles are: (i) the presence of shocks induces poor generalization and unstable gradients; (ii) entropy admissibility is not naturally enforced in the standard PINN framework; and (iii) boundary conditions in the hyperbolic regime are subtle and often imposed in physically inconsistent ways~\cite{fuks2021limitations, mishra2022physics}.

\medskip

\textbf{Contributions.} In this work, we introduce a new PINN-based framework designed to overcome these challenges. Our method is specifically tailored to solve initial-boundary value problems for nonlinear hyperbolic conservation laws in both scalar and system settings. The key ingredients of our approach are:

\begin{itemize}
  \item A \emph{parabolic regularization strategy}, introducing a small viscosity $\varepsilon > 0$ into the equation, allowing the network to learn a smooth solution $u^\varepsilon$ that converges, as $\varepsilon \to 0$, to the unique entropy solution of the original hyperbolic problem. This enables stable training and mitigates gradient pathologies caused by discontinuities.

  \item A \emph{weak imposition of boundary conditions} inspired by the Dubois–LeFloch theory~\cite{bib2}, which uses a nonlinear projection derived from the solution of local Riemann problems to enforce entropy consistency at the boundary. This replaces the usual Dirichlet loss term with a physically meaningful constraint.

  \item A \emph{general and modular framework} applicable not only to scalar equations but also to systems of conservation laws.
\end{itemize}

To the best of our knowledge, this is the first PINN architecture that combines viscosity-based regularization and entropy-consistent boundary enforcement within a unified and minimal neural framework. Despite its simplicity, based on standard feedforward networks and classical loss terms, the method achieves high accuracy in approximating entropy solutions, including shocks, contact discontinuities and boundary layer structures.

\medskip

\textbf{Numerical validation.} We demonstrate the effectiveness of our approach through detailed experiments on the inviscid Burgers equation with time-dependent inflow data, a canonical test case for nonlinear wave propagation. Reference solutions are computed using Godunov-type finite volume schemes equipped with the Dubois–LeFloch boundary projection, allowing precise quantitative comparisons. Results show that our method accurately captures the nonlinear wave dynamics and boundary layers, using shallow neural networks without any sophisticated regularization technique.

\medskip

\textbf{Structure of the paper.}  
The paper is organized as follows. In Section~\ref{sec1}, we introduce the mathematical setting and our approximation strategy, including the parabolic regularization and weak boundary formulation.
In Section~\ref{disc}, we compare our approach with the one developed in \cite{bib3}.
Then, Section~\ref{sec2} presents numerical experiments comparing our new PINN approach to a classical finite volume solver. Finally, Section~\ref{sec3} summarizes the findings and outlines directions for future research.

\section{Model Problem and Approximation Strategy}\label{sec1}

We aim to approximate the solution 
$\mathbf {u}: (x,t) \in  (a,b) \times \mathbb R_+^*
  \mapsto \mathbf {u} (x,t) \in \mathbb R ^p,\; p\geq 1,$ of the nonlinear strictly hyperbolic system 
\begin{equation}
\label{eq1}
 \partial _t \mathbf {u} + \partial_x \mathbf{f} (\mathbf {u}) 
= 0, \quad x \in (a,b), \quad t>0,
\end{equation}

with initial condition 
\begin{equation}
\label{init}
\mathbf{u} (x,0)=\mathbf{u} _0 (x), \quad x \in (a,b),
\end{equation}
and boundary conditions
\begin{equation}
\label{b1}
\mathbf{u} (a^+,t) = \mathbf{l}(t), \quad \mathbf{u} (b^-,t) = \mathbf{r}(t), \quad t> 0,
\end{equation}
where the initial data is assumed to have small total variation.

Recall that the nonconservative form of the nonlinear strictly hyperbolic system \eqref{eq1} is
\begin{equation}
\label{eq2}
 \partial _t \mathbf {u} + A(\mathbf {u}) \partial_x \mathbf {u} =0,
\end{equation}
where the matrix $A$ is the Jacobian of the (smooth) flux $\mathbf{f}$.
Now, let us consider the following parabolic regularization of~\eqref{eq2}:
\begin{equation}
\label{eq3}
 \partial _t \mathbf {u}^{\epsilon} + A(\mathbf {u}^{\epsilon}) \partial_x \mathbf {u}^{\epsilon} = \epsilon \partial_{x^2}^2 \mathbf {u}^{\epsilon}.
\end{equation}
The solutions of the Cauchy problem ~\eqref{eq3}--\eqref{init} are defined globally in time, satisfy uniform \( BV \) estimates that are independent of \( \epsilon \) and depend continuously on the initial data ~\eqref{init} in the \( L^1 \) norm, with a Lipschitz constant independent of \( t \) and \(\epsilon \). Moreover, letting \( \epsilon \to 0^+ \), these smooth solutions converge strongly in \( L^1_{loc} \), to the unique entropy weak solution of the original Cauchy problem~\eqref{eq1}--\eqref{init} (see, e.g.,~\cite{bib1}).

The notion of weak boundary conditions we adopt follows the framework proposed by Dubois and LeFloch \cite{bib2}, based on solutions to Riemann problems. More precisely, we consider the boundary conditions \eqref{b1} in the following weak sense:
\begin{subequations}
\label{new}
\begin{align}
\mathbf{u} (a^+,t) & \in \{\boldsymbol{\mathcal{W}} \left(0^+; \mathbf{l}(t) ; \mathbf{u}_R \right), \mathbf{u}_R \text{ varying in }  \mathbb R ^p \}, \quad t>0, \\
\mathbf{u} (b^-,t) & \in \{\boldsymbol{\mathcal{W}} \left(0^-;\mathbf{u}_L; \mathbf{r}(t) \right)  , \mathbf{u}_L \text{ varying in }  \mathbb R ^p \}, \quad t>0,
\end{align}
\end{subequations}
where $\boldsymbol{\mathcal{W}} ( \cdot; \cdot ; \cdot )$ denotes the self-similar solution of the Riemann problem associated to the nonlinear hyperbolic system \eqref{eq1} with the two last arguments as left and right initial states.

Finally, in practice, we approximate the entropy solution of the initial-boundary value problem~\eqref{eq1}--\eqref{init}--\eqref{b1} by solving the parabolic problem~\eqref{eq3}--\eqref{init}--\eqref{new} for a small viscosity parameter $\epsilon > 0$, using the classical physics-informed neural network approach \cite{bib4}.

\section{Discussion}\label{disc}

The proposed framework is not only faster and simpler, but is also more general than the one introduced in \cite{bib3}.

In particular, the algorithm developed by Chaumet and Giesselmann \cite{bib3} involves the approximate solution of a Min-Max optimization problem. As a result, the training procedure is significantly more computationally demanding than in our approach, which requires only the solution of a standard minimization problem.

With respect to the numerical treatment of boundary conditions, our method explicitly handles both scalar and system cases, in contrast to the algorithm in \cite{bib3}, which considers only the scalar setting.
However, the framework developed in \cite{bib3} can be extended to systems with boundary conditions by adopting the Riemann-problem-based formulation introduced in this work.

\section{Numerical Experiments}\label{sec2}

Our goal is to demonstrate that the proposed PINN framework, despite its simplicity, is capable of approximating physically consistent entropy solutions to nonlinear conservation laws with boundary interactions.
We thus assess the ability of the proposed method to approximate entropy solutions of nonlinear scalar conservation laws with boundary layers.

In this section, we illustrate the ability of our proposed algorithm to approximate the physically relevant entropy solution of the nonlinear hyperbolic Burgers equation, posed as a mixed initial-boundary value problem. Specifically, we aim to solve
\[
\partial_t u +  \partial_x \frac{u^2}{2} = 0,
\]
with given initial and boundary conditions, by learning the solution of its parabolic regularization:
\[
\partial_t u + u \partial_x u = \epsilon \partial_{xx} u,
\]
with \(\epsilon = 0.01\), on the space-time domain \((x,t) \in [-1,1] \times [0,1]\). The left boundary condition is time-dependent:
\[
u(-1^+, t) = t - 0.5,
\]
and is enforced in the weak sense. The right boundary at \(x = +1^-\) does not require explicit enforcement.

The goal is to recover the correct hyperbolic behavior, including shock formation, rarefaction waves, and nonlinear steepening, while enforcing boundary conditions in a weak form inspired by the Dubois–LeFloch theory~\cite{bib2}.

\subsection{Neural Network Architecture and Loss Function}

We employ a fully connected feedforward neural network with three hidden layers of 20 neurons each and hyperbolic tangent activation functions. The input is a space-time pair \((x,t)\), and the output is a scalar prediction \(u_\theta(x,t)\) parameterized by the network weights \(\theta\).

The total loss minimized during training is defined as:
\[
\mathcal{L} = \mathcal{L}_{\mathrm{PDE}} + \mathcal{L}_{\mathrm{IC}} + \mathcal{L}_{\mathrm{BC}}^{\mathrm{weak}},
\]
with the following components:
\begin{align*}
\mathcal{L}_{\mathrm{PDE}} &= \frac{1}{N_f} \sum_{i=1}^{N_f} \left( \partial_t u_\theta(x_i,t_i) + u_\theta(x_i,t_i) \, \partial_x u_\theta(x_i,t_i) \right. \\
&\left. \qquad\qquad\quad - \epsilon \, \partial^2_{x^2} u_\theta(x_i,t_i) \right)^2, \\
\mathcal{L}_{\mathrm{IC}} &= \frac{1}{N_0} \sum_{j=1}^{N_0} \left( u_\theta(x_j, 0) - u_0(x_j) \right)^2, \\
\mathcal{L}_{\mathrm{BC}}^{\mathrm{weak}} &= \frac{1}{N_{\mathrm{bc}}} \sum_{k=1}^{N_{\mathrm{bc}}} \left( u_\theta(-1, t_k) - \mathcal{W}(0^+; b(t_k), u_\theta(-1^+, t_k)) \right)^2.
\end{align*}

where \(b(t_k) = t_k - 0.5\) denotes the time-dependent inflow boundary data. Recall that \(\mathcal{W}(0^+; b, u_R)\) denotes the trace at \(x = 0^+\) of the entropy solution to the Riemann problem:
\[
\begin{cases}
\partial_t u + \partial_x \left(\frac{u^2}{2}\right) = 0, \\
u(x,0) = 
\begin{cases}
b, & x<0, \\
u_R, & x>0,
\end{cases}
\end{cases}
\]
which selects the unique entropy-admissible state according to the Dubois–LeFloch criterion~\cite{bib2}.

This weak enforcement allows the PINN to impose the correct physical state at the inflow, even in the presence of nonlinear interactions and shocks. The right boundary at \(x=1\) is left unconstrained, as all test cases are designed with right-going waves and no incoming characteristics at that boundary.

The network is trained using the Adam optimizer with a learning rate of 0.001 for 5000 epochs. We use \(N_f = 2000\) collocation points uniformly sampled in the interior of the space-time domain \([-1,1] \times [0,1]\), \(N_0 = 100\) initial condition points placed on a uniform grid over \([-1,1]\), and \(N_{\mathrm{bc}} = 100\) time points uniformly drawn from \([0,1]\) to impose the left boundary condition at \(x = -1\). No data is enforced at the right boundary \(x=1\), which acts as a natural outflow.

The training points \((x_i, t_i)\), \((x_j, 0)\), and \((-1, t_k)\) are drawn uniformly from the respective subdomains: the interior of \([-1,1]\times[0,1]\), the initial slice \(t=0\), and the left boundary. No data is used at the right boundary, consistent with the outgoing nature of the solution.

\subsection{Test Cases}

We consider three classical initial-boundary configurations that illustrate key features of nonlinear wave propagation:

\begin{itemize}
    \item \textbf{Test case 1 (shock formation):}
    \[
    u_0(x) = 
    \begin{cases}
    1, & x < 0, \\
    0, & x \geq 0,
    \end{cases}
    \quad u(1^-, t) = 0.
    \]
    A right-going shock emerges due to the convergence of characteristics and the jump in initial data.

    \item \textbf{Test case 2 (shock formation from smooth data):}
    \[
    u_0(x) = -\sin(\pi x), \quad u(1^-, t) = 0.
    \]

Although the initial data is smooth, the non-monotonic profile induces a compression of characteristics near \( x = 0 \), leading to the formation of a shock in finite time. Indeed, according to the method of characteristics, a discontinuity emerges around time \( t \approx 1/\pi \). This test case illustrates a fundamental phenomenon in nonlinear hyperbolic conservation laws: the spontaneous formation of discontinuities from smooth initial profiles.

    \item \textbf{Test case 3 (rarefaction wave):}
    \[
    u_0(x) = 
    \begin{cases}
    -1, & x < 0, \\
    1, & x \geq 0,
    \end{cases}
    \quad u(1^-, t) = 1.
    \]
    The discontinuity gives rise to a rarefaction fan with diverging characteristics and a smooth intermediate profile.
\end{itemize}

In each test case, we refrain from reporting explicit prediction errors, and instead assess the performance of the trained PINN by directly comparing its output (denoted as the numerical solution) to a high-resolution reference solution. The latter is computed using the finite volume method introduced by Dubois and LeFloch~\cite{bib2}, which combines the Godunov flux at internal cell interfaces with an entropy-consistent boundary condition at the inflow. More precisely, the reference solution is obtained on a uniform spatial grid of 5000 points over the interval \([-1,1]\), with a time step \(\Delta t = 0.01\, \Delta x\), up to the final time \(t = 0.75\). The numerical flux is consistent with the inviscid Burgers flux \(f(u) = \frac{u^2}{2}\), and the left boundary condition is imposed via a Riemann-based entropy projection as prescribed in~\cite{bib2}, while the right boundary is left free. This high-fidelity solver provides a robust and physically accurate benchmark against which we evaluate the PINN's ability to capture key nonlinear structures, including shocks, rarefaction waves, and boundary layers.

We now compare the PINN predictions to reference solutions computed via a high-resolution finite volume scheme, described in detail below. This comparison allows us to assess the ability of the network to capture nonlinear features such as shocks, rarefaction fans, and steep gradients, without relying on explicit knowledge of characteristic speeds or entropy formulations.

\subsection{Results and Physical Insights}

\begin{figure}[h]
    \centering
    \begin{subfigure}[t]{0.48\linewidth}
        \centering
        \includegraphics[width=\linewidth]{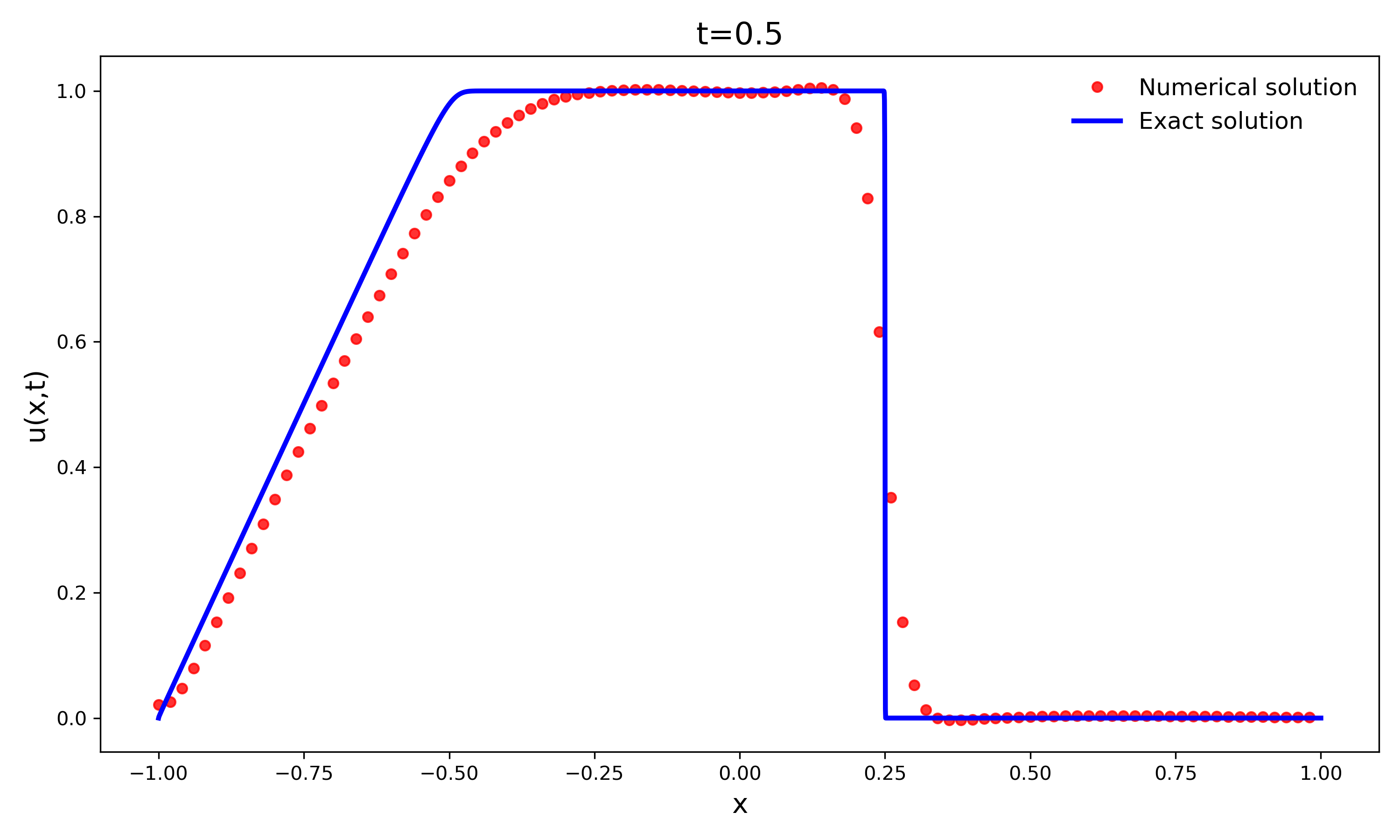}
        \caption{Solution to test case 1 at $t=0.5$}
    \end{subfigure}
    \hfill
    \begin{subfigure}[t]{0.48\linewidth}
        \centering
        \includegraphics[width=\linewidth]{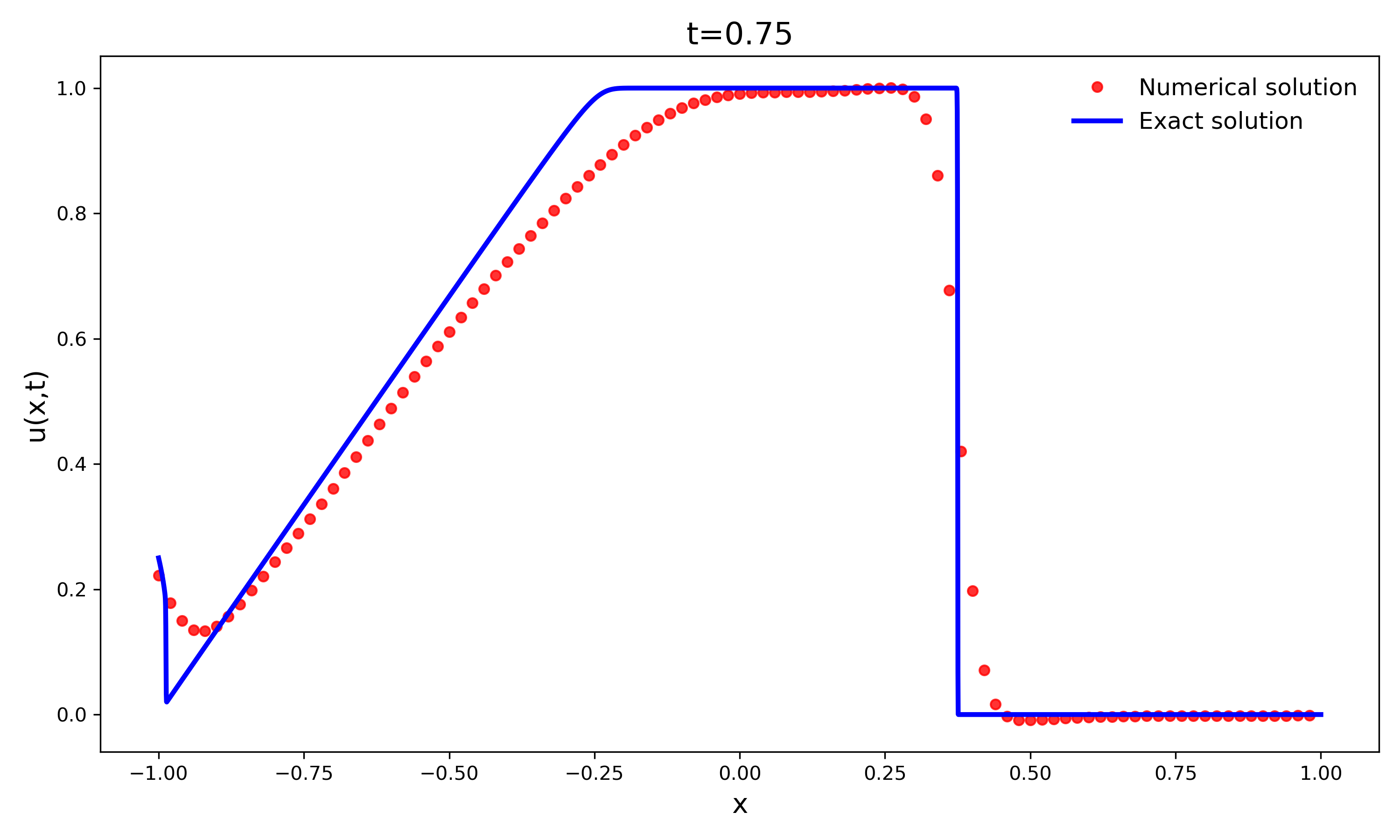}
        \caption{Solution to test case 1 at $t=0.75$}
    \end{subfigure}
    
    \vspace{0.4cm}
    
    \begin{subfigure}[t]{0.48\linewidth}
        \centering
        \includegraphics[width=\linewidth]{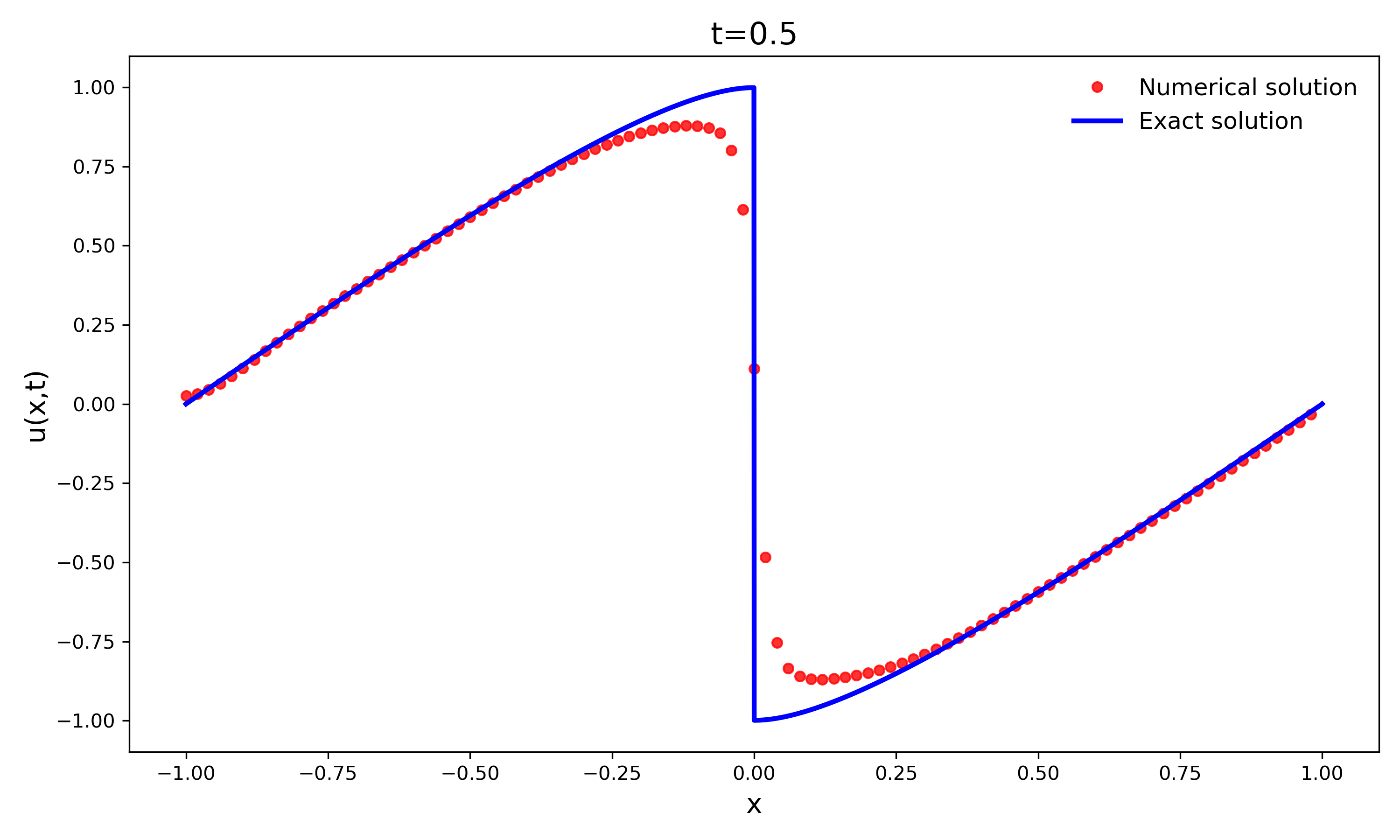}
        \caption{Solution to test case 2 at $t=0.5$}
    \end{subfigure}
    \hfill
    \begin{subfigure}[t]{0.48\linewidth}
        \centering
        \includegraphics[width=\linewidth]{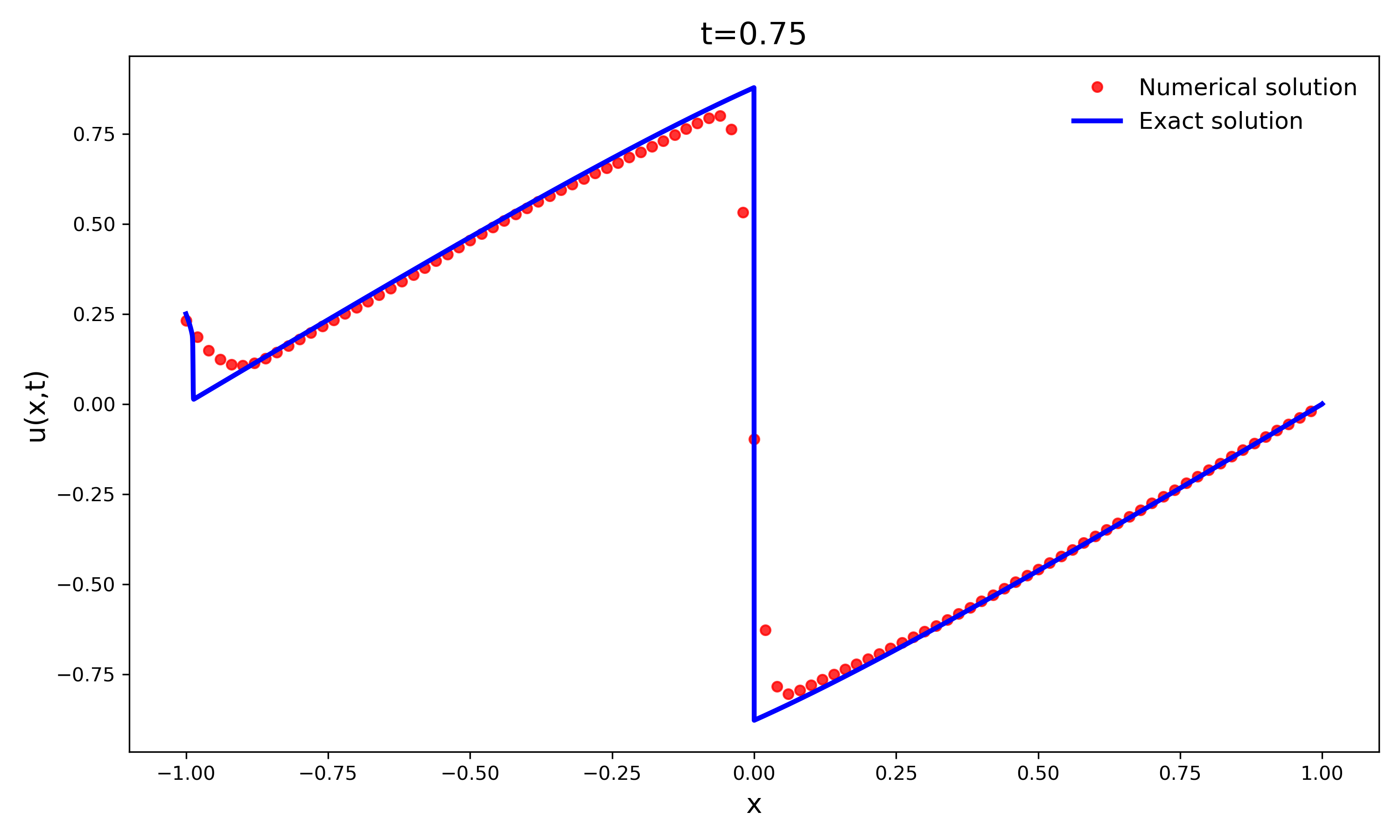}
        \caption{Solution to test case 2 at $t=0.75$}
    \end{subfigure}
    
    \begin{subfigure}[t]{0.48\linewidth}
        \centering
        \includegraphics[width=\linewidth]{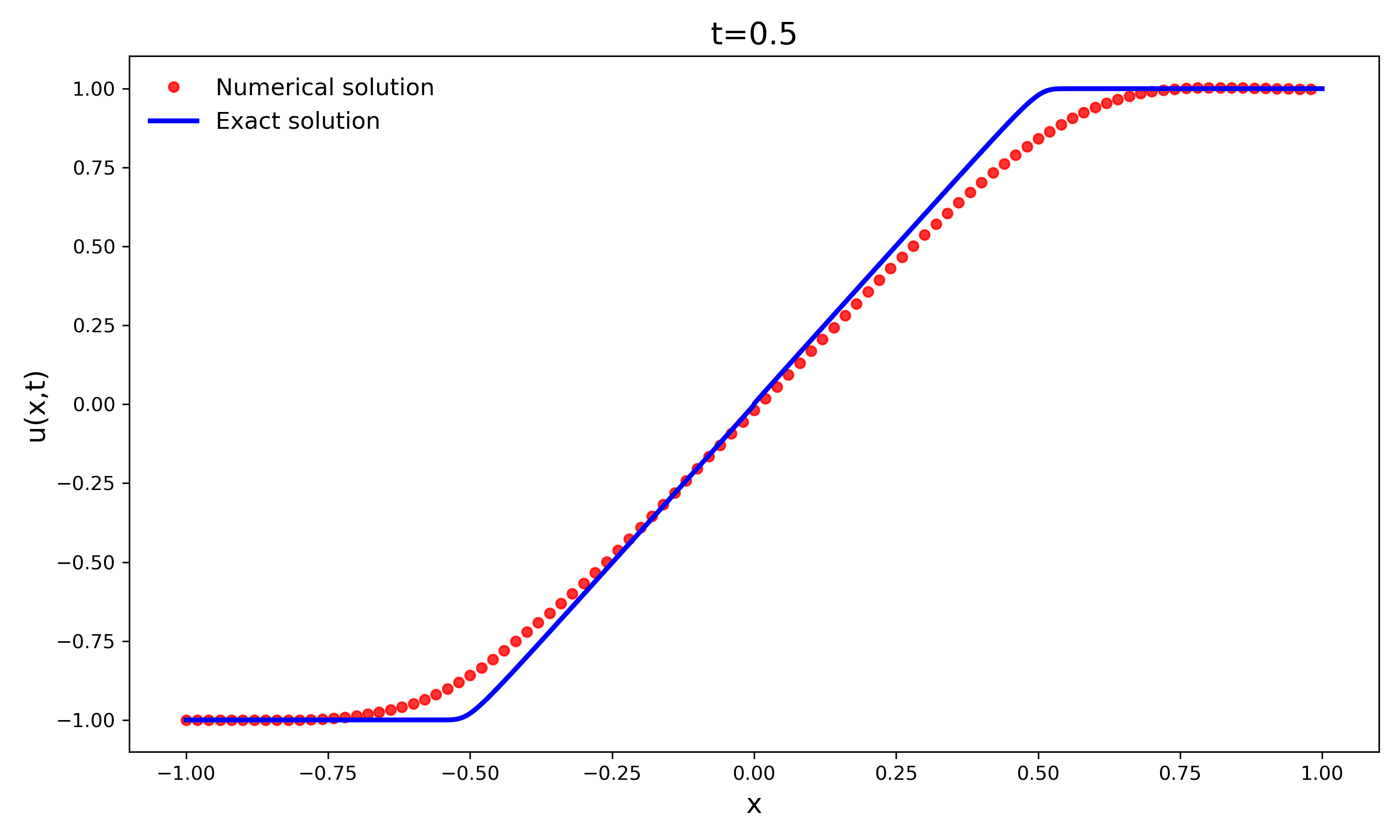}
        \caption{Solution to test case 3 at $t=0.5$}
    \end{subfigure}
    \hfill
    \begin{subfigure}[t]{0.48\linewidth}
        \centering
        \includegraphics[width=\linewidth]{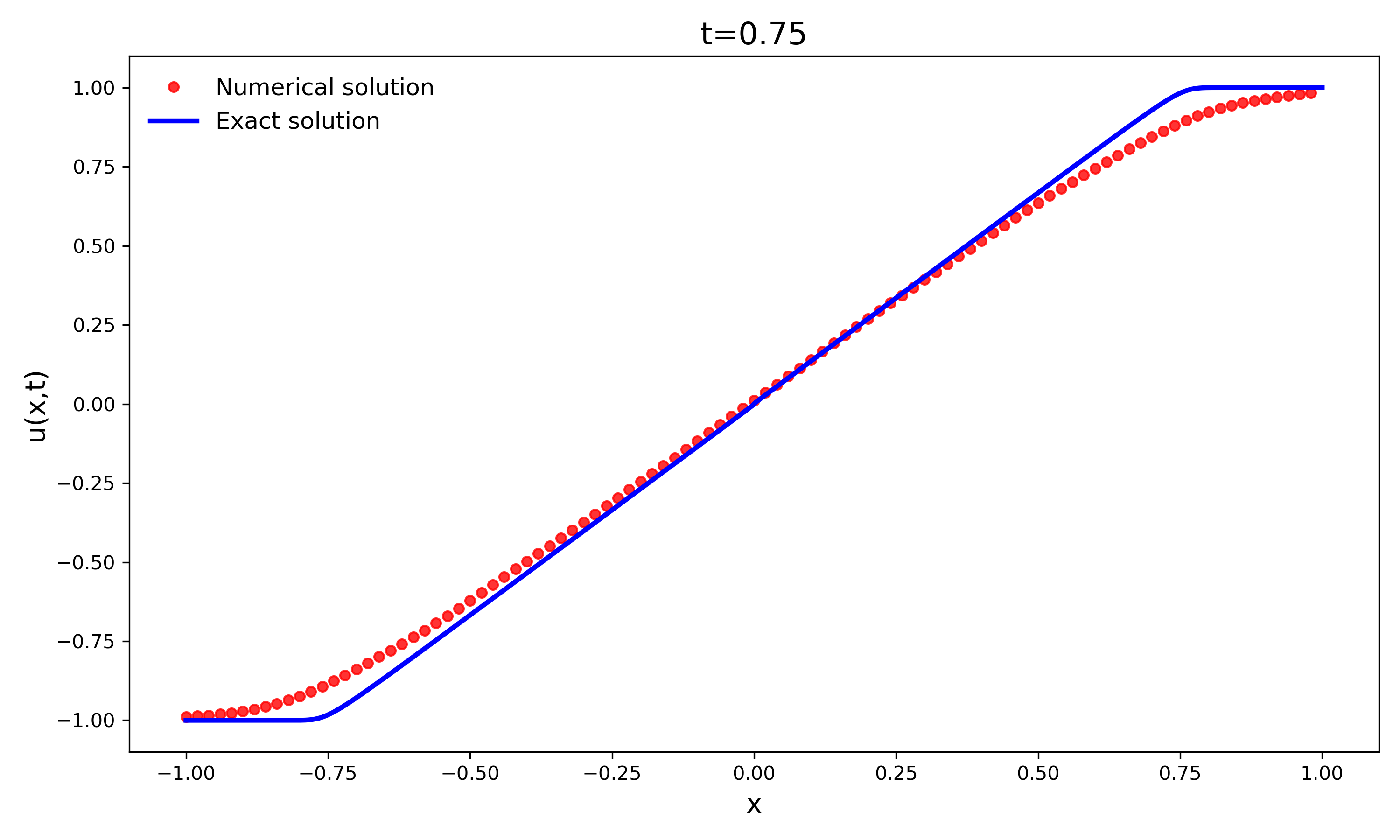}
        \caption{Solution to test case 3 at $t=0.75$}
    \end{subfigure}

    \caption{Comparison between neural network prediction and exact solution for each test case}
    \label{fig:resultats-2x2}
\end{figure}

Figure~\ref{fig:resultats-2x2} displays the predicted PINN solutions at times \(t = 0.5\) and \(t = 0.75\), compared to reference solutions.

\begin{itemize}
    \item \textbf{Test 1 (shock dynamics):}  
    The PINN accurately captures the formation and propagation of the shock wave. The left boundary acts as an active inflow, and the neural network dynamically adjusts the internal state to absorb and transmit the incoming signal. The use of the weak boundary projection \(\mathcal{W}\) ensures consistency with the entropy solution and prevents nonphysical oscillations.
    
    \item \textbf{Test 2 (shock formation from smooth data):}  
The network successfully captures the progressive steepening of the profile and the subsequent shock formation. The solution remains stable and accurate despite the presence of sharp gradients and emerging discontinuities. Notably, the PINN respects the entropy structure of the problem and avoids nonphysical oscillations near the shock, demonstrating its ability to handle singularity formation from smooth inputs.

    \item \textbf{Test 3 (rarefaction wave):}  
    The solution evolves into a smooth rarefaction fan. The network learns the self-similar expansion dynamics and captures the correct spreading speed. The result demonstrates the PINN’s ability to represent continuous but non-smooth solutions, driven purely by the PDE physics and the entropy condition.
\end{itemize}

These results illustrate the effectiveness of the proposed framework in capturing a wide range of nonlinear phenomena: boundary-driven wave injection, shock formation, and rarefaction development. The use of weak boundary enforcement based on local Riemann problems ensures physical admissibility, while the simple architecture and training protocol highlight the practicality and robustness of the method.

The full implementation of the neural solver and the reference finite volume scheme is provided in Appendix~\ref{appendix:code}.

These results suggest that physics-informed neural networks, when combined with entropy-consistent boundary treatments, can serve as robust solvers for nonlinear hyperbolic PDEs. 

\section{Concluding remarks}\label{sec3}

In this work, we have developed and validated a neural network-based framework for approximating entropy solutions to initial-boundary value problems for nonlinear strictly hyperbolic conservation laws. By leveraging parabolic regularization and a weak imposition of boundary conditions via local Riemann problems, our approach ensures physical consistency while maintaining computational simplicity. Numerical experiments based on the inviscid Burgers equation demonstrate that deep feedforward neural networks can accurately capture nonlinear wave dynamics and boundary layers, using minimal architectural complexity and no advanced regularization.

While the method is applicable to general systems in theory, future numerical investigations will aim to test its performance on systems of hyperbolic conservation laws. This framework also lends itself naturally to further generalizations. A natural extension is to incorporate nonlinear interface conditions and coupling constraints between adjacent domains. Such developments are particularly relevant for simulating multiscale flows or heterogeneous media, where flux discontinuities or nonlinear jump conditions arise.

These directions are closely related to the coupling strategies investigated in~\cite{bib5} and references therein, where well-posedness and entropy stability were addressed through analytical and numerical tools. Embedding such coupling mechanisms into a learning-based framework may open the door to hybrid solvers that combine the interpretability and structure-preserving properties of classical schemes with the adaptability of neural networks. Ultimately, this line of research may contribute to the design of robust neural–numerical methods for PDEs with complex boundary conditions, interface dynamics, and multiscale features.

\appendix
\section*{Appendix A: Python IMPLEMENTATION}
\label{appendix:code}

\begin{lstlisting}[language=Python, basicstyle=\footnotesize\ttfamily, breaklines=true]
import numpy as np
import tensorflow as tf
from tensorflow import keras
import matplotlib.pyplot as plt

# Function to solve viscous Burgers' equation using a PINN with weak BC at x = -1
def solve_burgers_pinn():
    nu = 0.01
    model = keras.Sequential([
        keras.layers.Dense(20, activation='tanh', input_shape=(2,)),
        keras.layers.Dense(20, activation='tanh'),
        keras.layers.Dense(20, activation='tanh'),
        keras.layers.Dense(1, activation=None)
    ])

    def pinn_loss(model, X_f):
        X_f_tf = tf.convert_to_tensor(X_f, dtype=tf.float32)
        with tf.GradientTape(persistent=True) as tape2:
            tape2.watch(X_f_tf)
            with tf.GradientTape(persistent=True) as tape1:
                tape1.watch(X_f_tf)
                u_pred = model(X_f_tf)
            u_x = tape1.gradient(u_pred, X_f_tf)[:, 0:1]
            u_t = tape1.gradient(u_pred, X_f_tf)[:, 1:2]
        u_xx = tape2.gradient(u_x, X_f_tf)[:, 0:1]
        f = u_t + u_pred * u_x - nu * u_xx
        return tf.reduce_mean(tf.square(f))

    def initial_condition(x):
        return np.where(x < 0, -1.0, 1.0).reshape(-1, 1)

    x_f = np.random.uniform(-1., 1., (2000, 1))
    t_f = np.random.uniform(0., 1., (2000, 1))
    X_f = np.hstack((x_f, t_f))

    x_0 = np.linspace(-1., 1., 100).reshape(-1, 1)
    t_0 = np.zeros_like(x_0)
    X_0 = np.hstack((x_0, t_0))
    u_0 = initial_condition(x_0)

    t_bc = np.random.uniform(0, 1, (100, 1))
    X_bc_left = np.hstack((-np.ones_like(t_bc), t_bc))
    optimizer = keras.optimizers.Adam(learning_rate=0.001)
    b = t_bc - 0.5 * np.ones_like(t_bc)

    def compute_u_bc_left(model, t_bc, b):
        u1_pred = model(tf.convert_to_tensor(X_bc_left, dtype=tf.float32)).numpy()
        u_bd = np.zeros_like(t_bc)
        for i in range(len(t_bc)):
            u1_i = u1_pred[i, 0]
            b_i = b[i, 0]
            if b_i > 0:
                if u1_i > 0 or (u1_i <= 0 and u1_i**2 < b_i**2):
                    u_bd[i, 0] = b_i
                else:
                    u_bd[i, 0] = u1_i
            else:
                if u1_i >= 0:
                    u_bd[i, 0] = b_i
                else:
                    u_bd[i, 0] = u1_i
        return u_bd

    def train_step():
        with tf.GradientTape() as tape:
            u_bc_left = compute_u_bc_left(model, t_bc, b)
            loss_pde = pinn_loss(model, X_f)
            loss_ic = tf.reduce_mean(tf.square(model(tf.convert_to_tensor(X_0, dtype=tf.float32)) - u_0))
            loss_bc_left = tf.reduce_mean(tf.square(model(tf.convert_to_tensor(X_bc_left, dtype=tf.float32)) - u_bc_left))
            loss = loss_pde + loss_ic + loss_bc_left
        grads = tape.gradient(loss, model.trainable_variables)
        optimizer.apply_gradients(zip(grads, model.trainable_variables))
        return loss.numpy()

    for epoch in range(5000):
        train_step()

    x_test = np.linspace(-1., 1., 5000).reshape(-1, 1)
    t_test = 0.75 * np.ones_like(x_test)
    X_test = np.hstack((x_test, t_test))
    u_pred = model(tf.convert_to_tensor(X_test, dtype=tf.float32)).numpy()
    return u_pred

# Godunov finite volume scheme with Dubois-LeFloch entropy projection at x = -1
def solve_burgers_equation():
    def flux(u):
        return 0.5 * u**2

    def godunov(u0, dx, dt, Nt):
        u = u0.copy()
        for n in range(Nt):
            b = n*dt - 0.5
            if b > 0:
                if u[1] > 0 or (u[1] <= 0 and u[1]**2 < b**2):
                    u[0] = b
                else:
                    u[0] = u[1]
            else:
                if u[1] >= 0:
                    u[0] = 0
                else:
                    u[0] = u[1]
            uL = u[:-1]
            uR = u[1:]
            u_interface = np.where(uL > uR, 
                                   np.where((uL + uR)/2 > 0, uL, uR), 
                                   np.where(uL >= 0, uL, np.where(uR <= 0, uR, 0)))
            F_interface = flux(u_interface)
            u[1:-1] -= dt/dx * (F_interface[1:] - F_interface[:-1])
        return u

    x = np.linspace(-1, 1, 5000)
    dx = x[1] - x[0]
    dt = 0.01 * dx
    Nt = int(0.75 / dt)
    u0 = np.where(x < 0, -1., 1.)
    u_num = godunov(u0, dx, dt, Nt)
    return u_num

x = np.linspace(-1, 1, 5000)
u_num_pinn = solve_burgers_pinn()
u_exact = solve_burgers_equation()

plt.figure(figsize=(10, 6))
step = 50
plt.plot(x[::step], u_num_pinn[::step], 'o', color='red', label='Numerical solution', markersize=5, alpha=0.8)
plt.plot(x, u_exact, 'b-', label='Exact solution', linewidth=2.5)
plt.title('t = 0.75', fontsize=16)
plt.xlabel('x', fontsize=14)
plt.ylabel('u(x,t)', fontsize=14)
plt.legend(fontsize=12, loc='best', frameon=False, shadow=True)
plt.tight_layout()
plt.savefig("casTest3_t75.png", dpi=300)
plt.show()
\end{lstlisting}


%

\ifCLASSOPTIONcaptionsoff
  \newpage
\fi



%

\bibliographystyle{plain}
\bibliography{sample}
%




\end{document}